# Third and fourth order well-balanced schemes for the shallow water equations based on the CWENO reconstruction

M. J. Castro Díaz,[*] M. Semplice[†]

**Abstract**

High order finite volume schemes for conservation laws are very useful in applications, due to their ability to compute accurate solutions on quite coarse meshes and with very few restrictions on the kind of cells employed in the discretization. For balance laws, the ability to approximate up to machine precision relevant steady states allows the scheme to compute accurately, also on coarse meshes, small perturbations of such states.

In this paper we propose third and fourth order accurate finite volume schemes for the shallow water equations. The schemes have the well-balanced property thanks to a path-conservative approach applied to an appropriate non-conservative reformulation of the equations. High order accuracy is achieved by designing truly two-dimensional reconstruction procedures of the CWENO type.

The novel schemes are tested for accuracy, well-balancing and shown to maintain posivity of the water height on wet/dry transitions. Finally they are applied to simulate the Tohoku 2011 tsunami event.

## 1 Introduction

Finite volume schemes are widespread for the numerical approximation of conservation and balance laws, thus being very important in many fields of applications. In this approach, the task of the numerical scheme is to compute the cell averages of the conserved variables on a computational grid at time $t^n + \Delta t$, starting from their values at time $t^n$. This task ought to be accomplished keeping in mind different and sometimes counteracting aims: high order accuracy and (essentially) non oscillatory properties are the main ones, but physical soundness of the computed values and approximation at machine precision of steady-state solutions of balance laws are also prominent in applications.

---

[*]Universidad de Málaga, Spain. email: *castro@anamat.cie.uma.es*
[†]Università di Torino, Italy. email: *matteo.semplice@unito.it*

In particular, the shallow water equations (SWEs) are useful to model free surface gravity waves whose wavelength is much larger than the characteristic bottom depth. This is, for example, the case of tsunami waves. In such large scale problems, it is mandatory to take into account the curvature of the Earth. Usually, the Earth can be approximated by a sphere and the equations are written in spherical coordinates. Although the PDE system is similar to the SWEs in the plane using Cartesian coordinates, new source terms appear due to the change of variables. In this paper, we follow the formulation as well as the general form of the well-balanced scheme described earlier [COP17].

If the advantage of using a high order scheme on a coarser mesh rather than a first order one on a fine mesh is unquestionable, on the other hand, an high order time advancement scheme needs higher order spatial information at time $t^n$ that are not directly described by the cell averages. It is therefore necessary to employ a reconstruction procedure that can recover the required in-cell high order spatial information from the cell averages in the neighbourhood.

In this paper we focus on the derivation of a high order well-balanced finite volume solver for the shallow water equations in cartesian and spherical coordinates combining a new Central WENO (CWENO) fourth order reconstruction operator, and a first order well balanced path-conservative scheme.

The quest for multi-dimensional non-oscillatory reconstruction procedures has been long standing. The Weighted Essentially Non Oscillatory reconstruction (WENO) [JS96, Shu09] is targeted at providing high-order accurate yet non-oscillatory point values of the conserved variables at cell boundaries. WENO relies on the existence (and positivity) of a set of, so called, linear or optimal weights that define a convex combination of low-degree interpolants that reproduce, at a specific point, the value of a high-degree central interpolant. Its implementation is very efficient for conservation laws in one-dimension or multidimensional Cartesian meshes. The values of the linear weights in WENO are fixed by accuracy requirements and depend on the location of the reconstruction point and on the size and relative location of the neighbouring cells. This has been a problem for example in the extensions to unstructured meshes, where various authors had to choose between using only low order polynomials and very complex computations of the optimal linear weights for each reconstruction point [HS99, SHS02, ZS09] or to employ a central high order polynomial with low order directionally-biased ones [DKTT07, BRDM09, TTD11].

For the integration of source terms in schemes for balance laws, some reconstruction points are located in the interior of the cell. The existence and positivity of the WENO linear weights is not guaranteed in this general situation [QS02]. On the other hand, the older ENO approach, which simply selects a reconstruction polynomial among a set of candidates, does not depend on such quantities and additionally provides a polynomial that is defined and uniformly accurate in the whole cell, without the need to compute a new set of linear and nonlinear coefficients for each reconstruction point. However the stencil of an ENO reconstruction is much wider that the stencil of a WENO one of the same order.

The Central WENO (CWENO) reconstruction, first introduced in the one-dimensional context by Levy, Puppo and Russo [LPR99], enjoys the good points of both ENO and WENO. In fact, CWENO makes use of linear weights, but their value is not determined by accuracy requirements and thus they can be fixed independently of the reconstruction point. As a consequence, CWENO, unlike WENO, does not suffer from the existence and positivity issues of the weights and moreover, like ENO, it yields a reconstruction polynomial that is valid in the entire cell. The coefficients actually employed in the reconstruction (called nonlinear weights) are then derived from the optimal values by a non-linear procedure whose task is to discard any information that leads to oscillatory polynomials. This procedure is very similar to the WENO one, but needs to be applied only once per cell and not once per reconstruction point. This is particularly advantageous for very high order schemes for conservation laws, due to the high number of quadrature points for the flux computations, and even more for balance laws, due to the additional quadrature points inside the cell (see e.g. [BRS18]).

After the above mentioned paper, the one-dimensional CWENO technique was extended to fifth order [Cap08, Zah09], the properties of the third order versions were studied in detail on uniform meshes [Kol14] and non-uniform ones [CS16] and finally arbitrary high order variants were introduced [CPSV18b, CPSV18a].

The CWENO procedure in more than one space dimension has been first employed in the form of a combination of a central parabola in two variables and four linear polynomials [LPR00a]. A similar construction, but with different definition of the polynomials, was exploited to obtain a third-order accurate reconstruction on two-dimensional quadrangular meshes, locally refined in a non-conforming fashion in a quad-tree type grid [SCR16]. A different approach was explored in two [LPR00b] and three [LP12] space dimension by defining a reconstruction as a convex combination (with nonlinear weights) of the central interpolating polynomials defined in each cell of the stencil, giving rise to reconstructions with quite large stencils compared to their order of accuracy. Also other authors defined two-dimensional reconstructions by convex combination of four polynomials, each of which can be either quadratic or linear [GOdlAM11]. More recently, the original idea of combining high and low degree polynomials was exploited also to obtain an arbitrary high order CWENO construction for triangular and tetrahedral meshes has been introduced in the finite volume schemes as a seed for an ADER predictor [DBSR17] and later also exploited in the subcell limiter for a DG scheme [DBS18]. Another CWENO construction, in the dimensional splitting approach, has also been presented[ZQ17].

In this paper, in the context of two-dimensional Cartesian meshes, we revise, adapting it to the needs of shallow water simulations, the CWENO-2D reconstruction of order three[SCR16] and introduce two novel CWENO reconstructions of order four. In particular, §2 describes the CWENO reconstructions, §3 describe the models of shallow water flows that will be used as examples. Next, §4 describes the application of the CWENO reconstructions in the the numer-ical schemes for shallow water flows in cartesian and spherical coordinates, §5 presents various numerical tests and §6 summarizes the conclusions of the paper.

|     |     |  9  |     |     |
| --- | --- | --- | --- | --- |
|     |  1  |  2  |  3  |     |
| 10  |  4  |  0  |  5  | 11  |
|     |  6  |  7  |  8  |     |
|     |     | 12  |     |     |

Figure 1: Stencil for the CWENO reconstructions. $\mathcal{S}_0$ for CWENO3 is composed by the cells $\Omega_0, \ldots, \Omega_8$, while CWENO4 uses also $\Omega_9, \ldots, \Omega_{12}$.

## 2 CWENO

All CWENO reconstructions employed in this paper are based on very small central stencils $\mathcal{S}_i$ that include the finite volume $\Omega_i$ and some of its neighbours. On $\mathcal{S}_i$, we define a so-called *optimal polynomial* $P_{\mathsf{opt}}$ (of degree 2 and 3 respectively) that interpolates exactly the cell average $\overline{u}_i$ and, in a least-square sense, all the other cell averages in $\mathcal{S}_i$. In standard situations (away from domain boundaries and dry states), the stencil $\mathcal{S}_i$ is the $3 \times 3$ central stencil for CWENO3 and a diamond-shaped stencil of 13 cells for CWENO4 (see also Fig. 1).

The CWENO reconstructions of this paper also make use of four lower degree polynomials ($P_r$ for $r = 1, \ldots, 4$) that fit the cell averages in a sector $\mathcal{S}_{i,r}$ of the central stencil $\mathcal{S}_i$ that determines the optimal polynomial. Each substencil $\mathcal{S}_{i,r} \subset \mathcal{S}_i$ is biased in a particular spatial direction, so that at least one of these polynomials may have the chance of avoiding discontinuities in the data that can be present in $\mathcal{S}_i$. In particular, $P_1$ will fit the cell averages of the cell that intersect the first quadrant of a reference system aligned with the grid and with origin in the center of $\Omega_0$, namely $\Omega_2, \Omega_3, \Omega_5$ if $\deg(P_r) = 2$ and also $\Omega_9, \Omega_{11}$ when $\deg(P_r) = 2$. In general, $P_r$ will interpolate the cell averages in the $r$-th quadrant of the above-mentioned reference system. (See §2.1 and §2.2 for more details.)

The core of all CWENO reconstructions is the equation[LPR99]

$$P_{\mathsf{opt}}(x,y) = d_0 P_0(x,y) + \sum_{r=1}^{4} d_r P_r(x,y), \qquad (1)$$

that rewrites $P_{\text{opt}}$ as a convex combination, using the so-called *linear coefficients* $d_0, \ldots, d_4 \in (0,1)$ that we will fix later. Having fixed $P_{\text{opt}}$ and $P_1, \ldots, P_4$ by interpolation requirements, equation (1) can be interpreted as the definition of a polynomial $P_0(x,y)$, of the same degree as $P_{\text{opt}}$, namely:

$$P_0(x,y) := \frac{1}{d_0}\left[P_{\text{opt}}(x,y) - \sum_{r=1}^{4} d_r P_r(x,y)\right]. \qquad (2)$$

The CWENO reconstruction polynomial is then defined as

$$R(x,y) = \omega_0 P_0(x,y) + \sum_{r=1}^{4} \omega_r P_r(x,y), \qquad (3)$$

where the weights $\omega_r$ should be designed to give essentially non-oscillatory properties to the reconstruction procedure. In particular one would like that $\omega_r \approx d_r$ for $r = 0, \ldots, 4$ in the case of smooth data, but that $\omega_r \approx 0$ if the $r$-th polynomial is constructed from discontinous data.

In order to fullfill these requirements in this paper we employ the Jiang-Shu indicators[JS96], which can be extended to the present two-dimensional setting as follows. Let us denote the partial derivative $(\partial/\partial_x)^{\alpha_1}(\partial/\partial_y)^{\alpha_2}$ as $\partial_\alpha$, where we introduced the multi-index $\alpha = (\alpha_1, \alpha_2)$ and further denote $|\alpha| = \alpha_1 + \alpha_2$ To a polynomial $P$, we associate the indicator $I[P]$ defined by

$$I[P] := \sum_{|\alpha|=1}^{\deg(P)} h^{2|\alpha|-2} \iint_{\Omega_0} (\partial_\alpha P)^2 \, \mathrm{d}x\mathrm{d}y \qquad (4)$$

where $h$ is a representative length associated to the cell, for which we choose the diameter $h = \sqrt{\Delta x^2 + \Delta y^2}$. Note that the normalization has been chosen so that $I[P] = O(1)$ when the data interpolated by $P$ are discontinuous.

The indicators (4) have the property that $I[P] = \mathcal{O}(h^2)$ if $P$ interpolates smooth data. With the help of the oscillation indicators, the nonlinear weights are then defined as usual by

$$\alpha_r = \frac{d_r}{(I[P_r] + \epsilon(h))^2} \qquad \omega_r = \frac{\alpha_r}{\sum_{s=0}^{4} \alpha_s} \qquad \text{for } r = 0, \ldots, 4. \qquad (5)$$

Note that for equation (2) to make sense, it is only needed that $d_0 \neq 0$. However, for (5) to define positive nonlinear weights, one should always choose $d_r \in (0,1)$ such that $\sum_{r=0}^{4} d_r = 1$. Once this prescription is satisfied, then the accuracy of the reconstruction polynomial (3) is uniform across the whole cell and depends only on the regularity of the data in the neighbourhood of the cell. In particular, maximal accuracy is achieved for smooth data, without having to choose linear weights that depend on the reconstruction point. As a consequence, the polynomial $R$ of (3) can be computed without prior knowledge of the points in $\Omega_0$ where it will be evaluated. This is particularly important for high order schemes, since many reconstructions points will be employed to compute the fluxes across all cell interfaces and to compute the quadrature of the source term in each cell.

The presence of a positive $\epsilon$ at the denominator in the previous formula is essential in order to avoid a division by zero in the case of reconstruction of very flat data. However, in the one-dimensional case it was proven (for both WENO [ABBM11] and CWENO [Kol14, CS16]) that the value of $\epsilon$ can affect the convergence rate close to local extrema. Intuitively, when a local extrema is present in the stencil, some of the indicators will be of size $\mathcal{O}(h^4)$ and some others will be $\asymp h^2$: unless $\epsilon$ is at least as big as $h^2$, formula (5) will bias the nonlinear weights towards the lower order polynomials with smaller indicators, effectively leading to an order reduction. Albeit no precise study has been performed in the 2D case, it has already been reported [SCR16] that both the choices $\epsilon(h) = h$ and $\epsilon(h) = h^2$ yield lower reconstruction errors than using a constant value for $\epsilon$.

It is important to note that, at a difference from standard WENO reconstruction procedures, the choice of the linear coefficients here is not related to the accuracy of the reconstruction at any given point. This allows not only to define a single reconstruction polynomial $R(x,y)$ that provides the reconstruction in the entire cell, but also gives the freedom to adjust the linear coefficients to other special needs. This will be exploited in the treatment of the wet-dry transitions (see Section 4.1).

## 2.1 Reconstruction of order 3

The third order reconstruction employed in this paper has been presented before in a non-uniform grid context[SCR16], but we describe it here in order to give the explicit formulas for the coefficients in the Cartesian uniform grid case.

In this case the *optimal polynomial* is the degree 2 polynomial with cell average $\overline{u}_0$ in $\Omega_0$ and approximating in the least-squares sense the cell averages of the cells $\Omega_1, \ldots, \Omega_8$, see Fig. 1. It is convenient to use a basis for $\mathbb{P}_2(x,y)$ such that all non-constant basis elements have zero mean in the central cell $\Omega_0$:

$$\begin{aligned}
&\varphi_0(x,y) = 1 &&\varphi_3(x,y) = (x-x_0)^2 - \tfrac{\Delta x^2}{12} \\
&\varphi_1(x,y) = (x-x_0) &&\varphi_4(x,y) = (y-y_0)^2 - \tfrac{\Delta y^2}{12} \\
&\varphi_2(x,y) = (y-y_0) &&\varphi_5(x,y) = (x-x_0)(y-y_0).
\end{aligned}$$

Thus we can write

$$P_{\mathsf{opt}}(x,y) = \overline{u}_0 + \sum_{k=1}^{5} c_k^{(2)} \varphi_k(x,y)$$

and the equations expressing the interpolation properties for the cells in the neighbourhood are

$$\overline{u}_j = \langle P_{\mathsf{opt}} \rangle_{\Omega_j} = \overline{u}_0 + \sum_{k=1}^{5} c_k^{(2)} \langle \varphi_k \rangle_{\Omega_j}, \qquad \text{for } j=1,\ldots,8.$$

Here we have denoted by $\langle \cdot \rangle_{\Omega_j}$ the cell average of a quantity on $\Omega_j$. The unknown coefficients $c_k^{(2)}$ are thus the (unconstrained) least square solution of the linear system of equations

$$\sum_{k=1}^{5} c_k^{(2)} \langle \varphi_k \rangle_{\Omega_j} = \overline{u}_j - \overline{u}_0, \qquad \text{for } j = 1, \ldots, 8.$$

In a regular Cartesian mesh one may form the $8 \times 5$ matrix with elements $A_{kj} = \langle \varphi_k \rangle_{\Omega_j}$ and pre-compute $B = (A^T A)^{-1} A^T$, so that the set of $c^{(2)}$ coefficients can be obtained pre-multiplying by $B$ the vector with elements $\overline{u}_j - \overline{u}_0$. The aforementioned computation can of course be performed symbolically to obtain the following explicit solution:

$$c_1^{(2)} = \tfrac{1}{6\Delta x} \left[ (\overline{u}_3 - \overline{u}_1) + (\overline{u}_5 - \overline{u}_4) + (\overline{u}_8 - \overline{u}_6) \right]$$
$$c_2^{(2)} = \tfrac{1}{6\Delta y} \left[ (\overline{u}_1 - \overline{u}_6) + (\overline{u}_2 - \overline{u}_7) + (\overline{u}_3 - \overline{u}_8) \right]$$
$$c_3^{(2)} = \tfrac{1}{10\Delta x^2} \left[ (\overline{u}_1 - 2\overline{u}_2 + \overline{u}_3) + 3(\overline{u}_4 - 2\overline{u}_0 + \overline{u}_5) + (\overline{u}_6 - 2\overline{u}_7 + \overline{u}_8) \right]$$
$$c_4^{(2)} = \tfrac{1}{10\Delta y^2} \left[ (\overline{u}_1 - 2\overline{u}_4 + \overline{u}_6) + 3(\overline{u}_2 - 2\overline{u}_0 + \overline{u}_7) + (\overline{u}_3 - 2\overline{u}_5 + \overline{u}_8) \right]$$
$$c_5^{(2)} = \tfrac{1}{4\Delta x \Delta y} \left[ (\overline{u}_3 - \overline{u}_1) - (\overline{u}_8 - \overline{u}_6) \right].$$

Note that each coefficient is a weighted average of the values of second order accurate discrete partial derivatives that can be computed in the stencil. In this respect, the $P_{\mathsf{opt}}$ employed here is different than the one of [LPR00b], where the least-square approach was not used.

The next step is to compute four first degree polynomials $P_r^{(1)}(x,y)$ that satisfy $\left\langle P_r^{(1)} \right\rangle_{\Omega_0} = \overline{u}_0$ and that fit, in the least-squares sense, a set of three cell averages in the neighbourhood. In particular, we will chose their stencils as follows (see Fig. 1):

| $P_r$ | $P_1^{(1)}$ | $P_2^{(1)}$ | $P_3^{(1)}$ | $P_4^{(1)}$ |
|---|---|---|---|---|
| $\mathcal{S}_{0,r}$ | $\Omega_2, \Omega_3, \Omega_5$ | $\Omega_5, \Omega_7, \Omega_8$ | $\Omega_4, \Omega_6, \Omega_7$ | $\Omega_1, \Omega_2, \Omega_4$ |

By writing

$$P_r^{(1)} = \overline{u}_0 + \sum_{k=1}^{2} c_{r,k} \varphi_k(x,y)$$

and proceeding as before, one gets

| | $c_{r,1}$ | $c_{r,2}$ |
|---|---|---|
| $P_1^{(1)}$ | $\frac{1}{3\Delta x}\left[(\overline{u}_3 - \overline{u}_2) + 2(\overline{u}_5 - \overline{u}_0)\right]$ | $\frac{1}{3\Delta y}\left[(\overline{u}_3 - \overline{u}_5) + 2(\overline{u}_2 - \overline{u}_0)\right]$ |
| $P_2^{(1)}$ | $\frac{1}{3\Delta x}\left[(\overline{u}_8 - \overline{u}_7) + 2(\overline{u}_5 - \overline{u}_0)\right]$ | $\frac{1}{3\Delta y}\left[(\overline{u}_5 - \overline{u}_8) + 2(\overline{u}_0 - \overline{u}_7)\right]$ |
| $P_3^{(1)}$ | $\frac{1}{3\Delta x}\left[(\overline{u}_7 - \overline{u}_6) + 2(\overline{u}_0 - \overline{u}_4)\right]$ | $\frac{1}{3\Delta y}\left[(\overline{u}_4 - \overline{u}_6) + 2(\overline{u}_0 - \overline{u}_7)\right]$ |
| $P_4^{(1)}$ | $\frac{1}{3\Delta x}\left[(\overline{u}_2 - \overline{u}_1) + 2(\overline{u}_0 - \overline{u}_4)\right]$ | $\frac{1}{3\Delta y}\left[(\overline{u}_1 - \overline{u}_4) + 2(\overline{u}_2 - \overline{u}_0)\right]$ |

The last polynomial appearing in the CWENO3 reconstruction is

$$P_0(x,y) = \frac{1}{d_0}\left[P_{\text{opt}} - \sum_{r=1}^{4} d_r P_r^{(1)}\right] = \overline{u}_0 + \sum_{k=1}^{4} c_{0,k}\varphi_k(x,y),$$

where

$$k = 1,2: \quad c_{0,k} = \frac{1}{d_0}\left(c_k^{(2)} - \sum_{r=1}^{4} d_r c_{r,k}^{(1)}\right)$$

$$k = 3,4,5: \quad c_{0,k} = \frac{1}{d_0}c_k^{(2)}.$$

## 2.2 Reconstructions of order 4

The reconstructions of order 4 employs as $P_{\text{opt}} = P^{(3)}$, the degree 3 polynomial in two variables that interpolates $\overline{u}_0$ exactly and the 12 cell averages in the diamond stencil depicted in Fig 1 in the sense of least squares. Note that this is very different from the tecnique used in [LPR02], which was later extended to three space dimensions in [LP12], that leads to a larger stencil.

We thus consider a basis for the 10 dimensional space $\mathbb{P}_3$ of polynomials of degree 3 in two variables that is composed by the functions $\varphi_0, \ldots, \varphi_5$ already introduced before, augmented by

$$\begin{array}{ll} \varphi_6 = (x-x_0)^3 & \varphi_8 = (x-x_0)(y-y_0)^2 \\ \varphi_7 = (y-y_0)^3 & \varphi_9 = (x-x_0)^2(y-y_0). \end{array}$$

Considering $P^{(3)} = \sum_{k=0}^{9} c_k^{(3)} \varphi_k(x,y)$, the coefficients $c_k^{(3)}$ are thus the least square solution of the linear system of equations

$$\sum_{k=1}^{9} c_k^{(3)} \langle \varphi_k \rangle_{\Omega_j} = \overline{u}_j - \overline{u}_0, \qquad \text{for } j = 1, \ldots, 12.$$

In a regular Cartesian mesh one may form the $12 \times 9$ matrix with elements $A_{kj} = \langle \varphi_k \rangle_{\Omega_j}$ and pre-compute $B = (A^T A)^{-1} A^T$, so that the set of $c$ coefficients can be obtained pre-multiplying by $B$ the vector with elements $\overline{u}_j - \overline{u}_0$. Explicitly

one finds that

$$c_1^{(3)} = \tfrac{1}{48\Delta x}\left[36(\overline{u}_5 - \overline{u}_5) - 5(\overline{u}_{11} - \overline{u}_{10}) - (\overline{u}_8 - \overline{u}_6) - (\overline{u}_3 - \overline{u}_1)\right]$$

$$c_2^{(3)} = \tfrac{1}{48\Delta x}\left[36(\overline{u}_2 - \overline{u}_7) - 5(\overline{u}_9 - \overline{u}_{12}) - (\overline{u}_1 - \overline{u}_6) - (\overline{u}_3 - \overline{u}_8)\right]$$

$$c_3^{(3)} = \tfrac{1}{714\Delta x^2}\left[\begin{array}{c} 76(\overline{u}_{11} - 2\overline{u}_0 + \overline{u}_{10}) + 19(\overline{u}_5 - 2\overline{u}_0 + \overline{u}_4) \\ +17(\overline{u}_3 - 2\overline{u}_2 + \overline{u}_1) + 17(\overline{u}_6 - 2\overline{u}_7 + \overline{u}_8) \\ +32(\overline{u}_2 - 2\overline{u}_0 + \overline{u}_7) - 8(\overline{u}_9 - 2\overline{u}_0 + \overline{u}_{12}) \end{array}\right]$$

$$c_4^{(3)} = \tfrac{1}{714\Delta y^2}\left[\begin{array}{c} 76(\overline{u}_9 - 2\overline{u}_0 + \overline{u}_{12}) + 19(\overline{u}_2 - 2\overline{u}_0 + \overline{u}_7) \\ +17(\overline{u}_1 - 2\overline{u}_4 + \overline{u}_6) + 17(\overline{u}_3 - 2\overline{u}_5 + \overline{u}_8) \\ +32(\overline{u}_4 - 2\overline{u}_0 + \overline{u}_5) - 8(\overline{u}_{10} - 2\overline{u}_0 + \overline{u}_{11}) \end{array}\right]$$

$$c_5^{(3)} = \tfrac{1}{4\Delta x \Delta y}\left[(\overline{u}_3 - \overline{u}_1) - (\overline{u}_8 - \overline{u}_6)\right]$$

$$c_6^{(3)} = \tfrac{1}{12\Delta x^3}\left[(\overline{u}_{11} - \overline{u}_{10}) - 2(\overline{u}_5 - \overline{u}_4)\right]$$

$$c_7^{(3)} = \tfrac{1}{12\Delta y^3}\left[(\overline{u}_9 - \overline{u}_{12}) - 2(\overline{u}_2 - \overline{u}_7)\right]$$

$$c_8^{(3)} = \tfrac{1}{4\Delta x \Delta y^2}\left[(\overline{u}_3 - \overline{u}_1) - 2(\overline{u}_5 - \overline{u}_4) + (\overline{u}_8 - \overline{u}_6)\right]$$

$$c_9^{(3)} = \tfrac{1}{4\Delta x^2 \Delta y}\left[(\overline{u}_1 - \overline{u}_6) - 2(\overline{u}_2 - \overline{u}_7) + (\overline{u}_3 - \overline{u}_8).\right].$$

The first reconstruction of order four that we propose is based on specializing (1) to

$$P_{\mathsf{opt}}(x,y) = P^{(3)}(x,y) = d_0 P_0(x,y) + \sum_{r=1}^{4} d_r P_r^{(1)}(x,y), \tag{6}$$

where the first degree polynomials coincide with those already employed in the third order accurate reconstruction described before. The reconstruction obtained from the decomposition (6) will be denoted by P3/P1 later on. It follows the approach whereby a high degree central polynomial is combined with first degree polynomials with very small stencils, in order to enhance the non-oscillatory properties of the reconstruction [DBSR17].

In this paper we propose also another fourth order reconstruction that will be denoted by P3/P2 and that is taylored to enhance the accuracy on smooth solutions. This latter makes use of four second degree polynomials $P_r^{(2)}(x,y)$ that satisfy $\left\langle P_r^{(2)}\right\rangle_{\Omega_0} = \overline{u}_0$ and that fit, in the least-squares sense, a set of five cell averages in the neighbourhood. In particular, we will chose their stencils as follows (see Fig. 1):

| $P_r$ | $P_1^{(2)}$ | $P_2^{(2)}$ | $P_3^{(2)}$ | $P_4^{(2)}$ |
|---|---|---|---|---|
| $\mathcal{S}_{0,r}$ | $\Omega_2, \Omega_3, \Omega_5$ $\Omega_9, \Omega_{11}$ | $\Omega_5, \Omega_7, \Omega_8$ $\Omega_{11}, \Omega_{12}$ | $\Omega_4, \Omega_6, \Omega_7$ $\Omega_{10}, \Omega_{12}$ | $\Omega_1, \Omega_2, \Omega_4$ $\Omega_9, \Omega_{10}$ |

The coefficients of $P_r^{(2)}$ are given in the following table:

|  | $P_1^{(2)}$ | $P_2^{(2)}$ |
|---|---|---|
| $c_{r,1}^{(2)}$ | $\frac{1}{2\Delta x}[4(\overline{u}_5 - \overline{u}_0) - (\overline{u}_{11} - \overline{u}_0)]$ | $\frac{1}{2\Delta x}[4(\overline{u}_5 - \overline{u}_0) - (\overline{u}_{11} - \overline{u}_0)]$ |
| $c_{r,2}^{(2)}$ | $\frac{1}{2\Delta y}[4(\overline{u}_2 - \overline{u}_0) - (\overline{u}_9 - \overline{u}_0)]$ | $\frac{1}{2\Delta y}[4(\overline{u}_{12} - \overline{u}_0) - (\overline{u}_7 - \overline{u}_0)]$ |
| $c_{r,3}^{(2)}$ | $\frac{1}{2\Delta x^2}[\overline{u}_{11} + u_0 - 2\overline{u}_5]$ | $\frac{1}{2\Delta x^2}[\overline{u}_{11} + u_0 - 2\overline{u}_5]$ |
| $c_{r,4}^{(2)}$ | $\frac{1}{2\Delta y^2}[\overline{u}_9 + u_0 - 2\overline{u}_2]$ | $\frac{1}{2\Delta y^2}[\overline{u}_{12} + u_0 - 2\overline{u}_7]$ |
| $c_{r,5}^{(2)}$ | $\frac{1}{\Delta x \Delta y}[(\overline{u}_3 - \overline{u}_2) - (\overline{u}_5 - \overline{u}_0)]$ | $\frac{1}{\Delta x \Delta y}[(\overline{u}_5 - \overline{u}_0) - (\overline{u}_8 - \overline{u}_7)]$ |

and

|  | $P_3^{(2)}$ | $P_4^{(2)}$ |
|---|---|---|
| $c_{r,1}^{(2)}$ | $\frac{1}{2\Delta x}[4(\overline{u}_0 - \overline{u}_4) - (\overline{u}_0 - \overline{u}_{10})]$ | $\frac{1}{2\Delta x}[4(\overline{u}_0 - \overline{u}_4) - (\overline{u}_0 - \overline{u}_{10})]$ |
| $c_{r,2}^{(2)}$ | $\frac{1}{2\Delta y}[4(\overline{u}_0 - \overline{u}_7) - (\overline{u}_0 - \overline{u}_{12})]$ | $\frac{1}{2\Delta y}[4(\overline{u}_2 - \overline{u}_0) - (\overline{u}_9 - \overline{u}_0)]$ |
| $c_{r,3}^{(2)}$ | $\frac{1}{2\Delta x^2}[\overline{u}_{10} + u_0 - 2\overline{u}_4]$ | $\frac{1}{2\Delta x^2}[\overline{u}_{10} + u_0 - 2\overline{u}_4]$ |
| $c_{r,4}^{(2)}$ | $\frac{1}{2\Delta y^2}[\overline{u}_{12} + u_0 - 2\overline{u}_7]$ | $\frac{1}{2\Delta y^2}[\overline{u}_9 + u_0 - 2\overline{u}_2]$ |
| $c_{r,5}^{(2)}$ | $\frac{1}{\Delta x \Delta y}[(\overline{u}_0 - \overline{u}_4) - (\overline{u}_7 - \overline{u}_6)]$ | $\frac{1}{\Delta x \Delta y}[(\overline{u}_9 - \overline{u}_0) - (\overline{u}_2 - \overline{u}_0)]$ |

The P3/P2 reconstruction is based on specializing (1) to

$$P_{\text{opt}}(x,y) = P^{(3)}(x,y) = d_0 P_0(x,y) + \sum_{r=1}^{4} d_r P_r^{(2)}(x,y). \qquad (7)$$

The reconstruction obtained from the decomposition (7) follows more closely the spirit of the original WENO [JS96] and CWENO [CPSV18b] reconstructions, in that it employs as lower degree polynomials the polynomials of maximal accuracy that can be determined by the data in some "directionally-biased" stencils. Both recontructions introduced in this section are designed to achieve fourth order accuracy on smooth data, but it is expected that P3/P2 should yield lower errors on regular solutions, but produce slightly larger spurious oscillations close to discontinuities than P3/P1.

**Boundary treatment** Close to a boundary, in order to reduce the number of ghost cells for which one must obtain extrapolated values, it is possible to obtain a 4-th order reconstruction with a smaller stencil. For example, close to the top-left boundary of the domain one may use only $\overline{u}_0, \ldots, \overline{u}_8, \overline{u}_{11}, \overline{u}_{12}$ to determine an optimal polynomial of degree 3, the standard polynomial of degree 2 in the south-west direction, and three first degree polynomials in the other directions. In order to obtain reconstructions that do not require ghost cells, one may investigate the adaption to the two-dimensional case of a technique that avoids altogether the use of ghost cells[NKS18].

## 2.3 Regularity indicators

In order to complete the description of the two-dimensional CWENO reconstructions , we report here explicit formulas for the oscillation indicators. They can be readily obtained applying the definition (4) to a polynomial of the form

$$P(x,y) = \overline{u_0} + \sum_{r=1}^{9} c_r \varphi_r(x,y).$$

Due to the symmetry of the basis functions $\varphi_r$, in the expression of the integrand, the terms containing odd powers of $(x-x_0)$ or $(y-y_0)$ do not contribute to the value of the indicator and thus do not appear in the formulae below.

For linear polynomials, only coefficients $c_1$ and $c_2$ are nonzero and the regularity indicator turns out to be

$$I[P^{(1)}] = A\left[c_1^2 + c_2^2\right],$$

where $A = \Delta x \Delta y$. For second order polynomials we have

$$I[P^{(2)}] = I[P^{(1)}] + A\left[\Delta x^2 \tfrac{c_3^2}{3} + \Delta y^2 \tfrac{c_4^2}{3} + (\Delta x^2 + \Delta y^2)(4c_3^2 + 4c_4^2 + \tfrac{13}{12}c_5^2)\right].$$

Finally, for a third order polynomial one has

$$I[P] = I[P^{(2)}] + A \begin{bmatrix} \frac{\Delta x^2}{2} c_1 c_6 + \frac{\Delta y^2}{2} c_2 c_7 + \frac{\Delta x^2 \Delta y^2}{24}(c_6 c_8 + c_7 c_9) \\ + \frac{1}{6}(\Delta x^2 c_2 c_9 + \Delta y^2 c_1 c_8) \\ + \frac{3}{80}(1043 \Delta x^4 + 2000 \Delta x^2 \Delta y^2 + 960 \Delta y^4) c_6^2 \\ + \frac{3}{80}(960 \Delta x^4 + 2000 \Delta x^2 \Delta y^2 + 1043 \Delta y^4) c_7^2 \\ + \frac{1}{720}(3120 \Delta x^4 + 6260 \Delta x^2 \Delta y^2 + 3129 \Delta y^4) c_8^2 \\ + \frac{1}{720}(3129 \Delta x^4 + 6260 \Delta x^2 \Delta y^2 + 3120 \Delta y^4) c_8^2 \end{bmatrix}.$$

## 3  The shallow water equations

In this paper we apply the novel CWENO reconstruction to simulations with the well known 2D nonlinear one-layer Shallow-water system, which, in cartesian coordinates, reads

$$\begin{cases} \partial_t h + \partial_x q_x + \partial_y q_y = 0, \\ \partial_t q_x + \partial_x \left(\dfrac{q_x^2}{h} + \dfrac{g}{2} h^2\right) + \partial_y \left(\dfrac{q_x q_y}{h}\right) = gh \partial_x H \\ \partial_t q_y + \partial_x \left(\dfrac{q_x q_y}{h}\right) + \partial_y \left(\dfrac{q_x^2}{h} + \dfrac{g}{2} h^2\right) = gh \partial_y H. \end{cases} \qquad (8)$$

In the previous system, $h(\mathbf{x},t)$, denotes the thickness of the water layer at point $\mathbf{x} \in \Omega \subset \mathbb{R}^2$ at time $t$, being $\Omega$ the horizontal projection of the 3D water body. $H(\mathbf{x})$ is the depth of the bottom at point $\mathbf{x}$ measured from a fixed level of reference. Let us also define the function $\eta(\mathbf{x},t) = h(\mathbf{x},t) - H(\mathbf{x})$ that corresponds to the free surface of the fluid. Let us denote by $\vec{q}(\mathbf{x},t) = (q_x(\mathbf{x},t), q_y(\mathbf{x},t))$ the mass-flow of the water layer at point $\mathbf{x}$ at time $t$. The mass-flow is related to the height-averaged velocity $\vec{u}(\mathbf{x},t)$ by means of the expression: $\vec{q}(\mathbf{x},t) = h(\mathbf{x},t)\, \vec{u}(\mathbf{x},t)$.

The previous system can be rewritten as

$$\begin{cases} \partial_t h + \partial_x q_x + \partial_y q_y = 0, \\ \partial_t q_x + \partial_x \left(\dfrac{q_x^2}{h}\right) + \partial_y \left(\dfrac{q_x q_y}{h}\right) + gh\partial_x \eta = 0 \\ \partial_t q_y + \partial_x \left(\dfrac{q_x q_y}{h}\right) + \partial_y \left(\dfrac{q_y^2}{h}\right) + gh\partial_y \eta = 0, \end{cases} \qquad (9)$$

or in a more compact form as

$$\partial_t w + \partial_x F_x(w) + \partial_y F_y(w) + T_x^p(w)\partial_x \eta + T_y^p(w)\partial_y \eta = 0, \qquad (10)$$

where

$$w = \begin{bmatrix} h \\ q_x \\ q_y \end{bmatrix} \qquad (11)$$

$$F_x(w) = \begin{bmatrix} q_x \\ \dfrac{q_x^2}{h} \\ \dfrac{q_x q_y}{h} \end{bmatrix} \qquad F_y(w) = \begin{bmatrix} q_y \\ \dfrac{q_x q_y}{h} \\ \dfrac{q_y^2}{h} \end{bmatrix} \qquad (12)$$

$$T_x^p(w) = \begin{bmatrix} 0 \\ gh \\ 0 \end{bmatrix} \qquad T_y^p(w) = \begin{bmatrix} 0 \\ 0 \\ gh \end{bmatrix}. \qquad (13)$$

Note that $\vec{F} = (F_x, F_y)$ is the convective flux and $\vec{T^p}(w) \cdot \nabla \eta$ is the pressure term where $\vec{T^p} = (T_x^p, T_y^p)$.

Finally, let us remark that the stationary solutions corresponding to water at rest are given by

$$u_x = 0, \quad u_y = 0, \quad \eta = \bar{\eta}, \qquad (14)$$

where $\bar{\eta}$ is a constant associated to the elevation of the undisturbed water.

A more interesting system is obtained when we consider the one-layer shallow-water system on the sphere. Here we use the formulation proposed in [COP17]:

$$\begin{cases} \partial_t h_\sigma + \dfrac{1}{R}\left(\partial_\theta \left(\dfrac{Q_\theta}{\cos(\varphi)}\right) + \partial_\varphi Q_\varphi\right) = 0, \\ \partial_t Q_\theta + \dfrac{1}{R}\partial_\theta \left(\dfrac{Q_\theta^2}{h_\sigma \cos(\varphi)}\right) + \dfrac{1}{R}\partial_\varphi \left(\dfrac{Q_\theta Q_\varphi}{h_\sigma}\right) - \dfrac{Q_\theta Q_\varphi}{R h_\sigma}\tan(\varphi) + \dfrac{gh_\sigma}{R\cos^2(\varphi)}\partial_\theta \eta_\sigma = 0, \\ \partial_t Q_\varphi + \dfrac{1}{R}\partial_\theta \left(\dfrac{Q_\varphi Q_\theta}{h_\sigma \cos(\varphi)}\right) + \dfrac{1}{R}\partial_\varphi \left(\dfrac{Q_\varphi^2}{h_\sigma}\right) + \left(\dfrac{Q_\theta^2}{R h_\sigma} + \dfrac{gh_\sigma \eta_\sigma}{R\cos(\varphi)}\right)\tan(\varphi) + \dfrac{gh_\sigma}{R\cos(\varphi)}\partial_\varphi \eta_\sigma = 0. \end{cases}$$
$$(15)$$

In the previos system, $R$ is the radius; $\mathbf{x} = (\theta, \varphi)$, the longitude and latitude; $g$, the gravity; $h$ is the thickness of the water layer and $h_\sigma = h\cos(\varphi)$; $H$ is the bottom depth, $H_\sigma = H\cos(\varphi)$ and $\eta_\sigma = h_\sigma - H_\sigma$; and finally, $u_\theta$, $u_\varphi$, are the longitudinal and latitudinal velocities averaged in the normal direction and $Q_\varphi = \cos(\varphi) q_\varphi$, $Q_\theta = \cos(\varphi) q_\theta$, with $q_\varphi = hu_\varphi$ and $q_\theta = hu_\theta$.

**Remark 1.** *Since $\cos(\varphi)$ is continuous, it can be easily checked that the Rankine-Hugoniot conditions are the usual ones. Moreover, if $H$ is assumed to be smooth, then the products $h_\sigma \partial_\theta \eta_\sigma$ and $h_\sigma \partial_\varphi \eta_\sigma$ are well defined.*

Finally, denoting $\sigma = \cos(\varphi)$ the system can be written as follows:

$$\begin{cases} \partial_t h_\sigma + \dfrac{1}{R}\left(\partial_\theta\left(\dfrac{Q_\theta}{\sigma}\right) + \partial_\varphi Q_\varphi\right) = 0, \\ \partial_t Q_\theta + \dfrac{1}{R}\partial_\theta\left(\dfrac{Q_\theta^2}{h_\sigma \sigma}\right) + \dfrac{1}{R}\partial_\varphi\left(\dfrac{Q_\theta Q_\varphi}{h_\sigma}\right) + \dfrac{Q_\theta Q_\varphi}{Rh_\sigma \sigma}\partial_\varphi \sigma + \dfrac{gh_\sigma}{R\sigma^2}\partial_\theta \eta_\sigma = 0, \\ \partial_t Q_\varphi + \dfrac{1}{R}\partial_\theta\left(\dfrac{Q_\varphi Q_\theta}{h_\sigma \sigma}\right) + \dfrac{1}{R}\partial_\varphi\left(\dfrac{Q_\varphi^2}{h_\sigma}\right) - \left(\dfrac{Q_\theta^2}{Rh_\sigma \sigma} + \dfrac{gh_\sigma \eta_\sigma}{R\sigma^2}\right)\partial_\varphi \sigma + \dfrac{gh_\sigma}{R\sigma}\partial_\varphi \eta_\sigma = 0. \end{cases}$$
(16)

As in the case of Cartesian coordinates, we consider the more compact form

$$\partial_t w + \dfrac{1}{R}\left(\partial_\theta F_\theta(W) + \partial_\varphi F_\varphi(w) + T^p_\theta(W)\partial_\theta \eta_\sigma + T^p_\varphi(W)\partial_\varphi \eta_\sigma + G_\varphi(W)\partial_\varphi \sigma\right) = 0$$
(17)

where

$$W = \begin{bmatrix} w \\ \sigma \end{bmatrix} = \begin{bmatrix} h_\sigma \\ Q_\theta \\ Q_\varphi \\ \sigma \end{bmatrix}$$
(18)

$$F_\theta(W) = \begin{bmatrix} \dfrac{Q_\theta}{\sigma} \\ \dfrac{Q_\theta^2}{h_\sigma \sigma} \\ \dfrac{Q_\varphi Q_\theta}{h_\sigma \sigma} \end{bmatrix} \quad F_\varphi(w) = \begin{bmatrix} Q_\varphi \\ \dfrac{Q_\theta Q_\varphi}{h_\sigma} \\ \dfrac{Q_\varphi^2}{h_\sigma} \end{bmatrix}$$
(19)

$$T^p_\theta(W) = \begin{bmatrix} 0 \\ \dfrac{gh_\sigma}{\sigma^2} \\ 0 \end{bmatrix} \quad T^p_\varphi(W) = \begin{bmatrix} 0 \\ 0 \\ \dfrac{gh_\sigma}{\sigma} \end{bmatrix}$$
(20)

$$G_\varphi(W) = G^1_\varphi(W) + G^2_\varphi(h_\sigma, \eta_\sigma, \sigma)$$
(21)

$$G^1_\varphi(W) = \begin{bmatrix} 0 \\ \dfrac{Q_\theta Q_\varphi}{h_\sigma \sigma} \\ -\dfrac{Q_\theta^2}{h_\sigma \sigma} \end{bmatrix} \quad G^2_\varphi(h_\sigma, \eta_\sigma, \sigma) = \begin{bmatrix} 0 \\ 0 \\ -\dfrac{gh_\sigma \eta_\sigma}{\sigma^2} \end{bmatrix}.$$
(22)

The advective flux $\vec{F} = (F_\theta, F_\varphi)$ and the pressure term $\vec{T}^p = (T_\theta^p, T_\varphi^p)$ satisfy the following rotational invariance-like properties: given a state $W$ and a unit vector $\vec{n} = [n_\theta, n_\varphi]^T$ one has

$$R_{\vec{\nu}} F_{\vec{n}}(W) = \delta\, F(R_{\vec{\nu}} w), \quad R_{\vec{\nu}} T_{\vec{n}}^p(W) = \delta\, T^p\left(\frac{h_\sigma}{\sigma}\right), \qquad (23)$$

where

$$F_{\vec{n}}(W) = n_\theta F_\theta(W) + n_\varphi F_\varphi(w), \quad T_{\vec{n}}^p(W) = n_\theta T_\theta^p(W) + n_\varphi T_\varphi^p(W), \qquad (24)$$

$$\delta = \sqrt{\frac{n_\theta^2}{\sigma^2} + n_\varphi^2}, \quad \vec{\nu} = \begin{bmatrix} \nu_\theta \\ \nu_\varphi \end{bmatrix} = \begin{bmatrix} \dfrac{n_\theta}{\sigma\delta} \\ \dfrac{n_\varphi}{\delta} \end{bmatrix}, \quad \vec{\nu}^\perp = \begin{bmatrix} -\nu_\varphi \\ \nu_\theta \end{bmatrix}, \qquad (25)$$

$$R_{\vec{\nu}} = \begin{bmatrix} 1 & 0 & 0 \\ 0 & \nu_\theta & \nu_\varphi \\ 0 & -\nu_\varphi & \nu_\theta \end{bmatrix}, \qquad (26)$$

and, for every $U_\sigma = [h_\sigma, Q_{\vec{\nu}}, Q_{\vec{\nu}^\perp}]^T$ and $h$:

$$F(U_\sigma) = \begin{bmatrix} Q_{\vec{\nu}} \\ \dfrac{Q_{\vec{\nu}}^2}{h_\sigma} \\ \dfrac{Q_{\vec{\nu}} Q_{\vec{\nu}^\perp}}{h_\sigma} \end{bmatrix}, \quad T^p(h) = \begin{bmatrix} 0 \\ gh \\ 0 \end{bmatrix}. \qquad (27)$$

Finally, let us remark that the stationary solutions corresponding to water at rest are given by

$$u_\varphi = 0, \quad u_\theta = 0, \quad \eta_\sigma = \bar{\eta} \cos(\varphi), \qquad (28)$$

where $\bar{\eta}$ is a constant associated to the elevation of the undisturbed water.

## 4 Numerical scheme

Here we briefly describe the discretization of the shallow-water system on the sphere proposed in [COP17]. Note that the shallow-water system in cartesian coordinates is a particular case where $\sigma = 1$ and $R = 1$. Therefore, the numerical scheme that follows can be easily adapted to that particular case.

Given a domain $\Omega$ in the $\theta$-$\varphi$ plane, we consider a partition in cells $\mathcal{T} = \{\Omega_i\}_{i=1}^{NC}$ with the usual properties: the cells are closed convex polygons and the intersection of two cells can only be a vertex, an edge, or empty. Two cells $\Omega_i$ and $\Omega_j$ are said to be neighbours it they share an edge $E_{i,j}$. The unit vector normal to $E_{i,j}$ pointing from $\Omega_i$ to $\Omega_j$ is denoted by $\vec{n}_{i,j} = [n_{i,j}^\theta, n_{i,j}^\varphi]^T$. Finally, $|\Omega_i|$ and $|E_{i,j}|$ represent the area of $\Omega_i$ and the length of $E_{i,j}$ respectively. Given a cell $\Omega_i$, we denote by $\mathcal{N}_i$ the set of indexes of the Neumann neighbours of $\Omega_i$. Finally, $\Delta \mathbf{x}$ represents the maximum of the diameters of the cells.

The approximation of the average of the solution at the cell $\Omega_i$ at time $t$ will be represented by:
$$w_i(t) = \begin{bmatrix} h_{\sigma,i}(t) \\ Q_{\theta,i}(t) \\ Q_{\varphi,i}(t) \end{bmatrix},$$

and
$$W_i(t) = \begin{bmatrix} w_i(t) \\ \sigma_i \end{bmatrix},$$

where $\sigma_i$ is the exact average of $\sigma = \cos(\varphi)$.

According to [COP17], a high order well-balanced numerical scheme for the water at rest solution (28) for system (10) is given by:

$$\begin{aligned} w_i'(t) &= -\frac{1}{R|\Omega_i|} \Bigg( \sum_{j \in \mathcal{N}_i} \left( \int_{E_{i,j}} F_{\vec{n}_{i,j}}(W_{i,j}^-(\gamma)) \, d\gamma + \int_{E_{i,j}} \delta_{i,j}(\gamma) D_{i,j}^-(\gamma) \, d\gamma \right) \\ &\quad + \int_{\Omega_i} \left( T_\theta^p(P_i(\mathbf{x})) \partial_\theta p_{f,i}(\mathbf{x}) + T_\varphi^p(P_i(\mathbf{x})) \partial_\varphi p_{f,i}(\mathbf{x}) \right) d\mathbf{x} \quad (29) \\ &\quad + \int_{\Omega_i} \left( G_\varphi^1(P_i(\mathbf{x})) + G_\varphi^2(p_{h_\sigma}(\mathbf{x}), p_{f,i}(\mathbf{x}), \sigma(\mathbf{x})) \right) \partial_\varphi \sigma(\mathbf{x}) \, d\mathbf{x} \Bigg), \end{aligned}$$

where the equality
$$T_\theta^p(P_i(\mathbf{x}))\partial_\theta(\bar{\eta}\cos(\varphi)) + T_\varphi^p(P_i(\mathbf{x}))\partial_\varphi(\bar{\eta}\cos(\varphi)) + G_\varphi^2(p_{h_\sigma}, \bar{\eta}\cos(\varphi), \sigma)\partial_\varphi \cos(\varphi) = 0 \quad (30)$$

has been used.

Notice that in the case of the shallow-water equation in cartesian coordinates the last integral vanishes as $\sigma$ is constant.

In the previous numerical scheme, the following notation has been used:

$$P_i(\mathbf{x}) = \begin{bmatrix} p_{h_\sigma,i}(\mathbf{x}) \\ p_{Q_\theta,i}(\mathbf{x}) \\ p_{Q_\varphi,i}(\mathbf{x}) \\ \sigma(\mathbf{x}) \end{bmatrix}, \quad W_{i,j}^\pm(\gamma) = \begin{bmatrix} w_{i,j}^\pm(\gamma) \\ \sigma(\gamma) \end{bmatrix} = \begin{bmatrix} h_{\sigma,i,j}^\pm(\gamma) \\ Q_{\theta,i,j}^\pm(\gamma) \\ Q_{\varphi,i,j}^\pm(\gamma) \\ \sigma(\gamma) \end{bmatrix}, \quad (31)$$

represent the reconstructions of $\{h_{\sigma,i}\}$, $\{Q_{\theta,i}\}$, $\{Q_{\varphi,i}\}$ at $\mathbf{x} \in \mathring{\Omega}_i$ and at $\gamma \in E_{i,j}$ respectively; $p_{\eta_\sigma,i}(\mathbf{x})$ and $\eta_{\sigma,i,j}^\pm(\gamma)$ represent the reconstruction of $\{\eta_{\sigma,i}\}$ at $\mathbf{x} \in \mathring{\Omega}_i$ and at $\gamma \in E_{i,j}$ respectively; $p_{f,i}(\mathbf{x})$ is a high order reconstruction of the fluctuations of the free surface with respect to the stationary solutions corresponding to water at rest (28), and it will be described later;

$$\delta_{i,j}(\gamma) = \sqrt{\left(\frac{n_{i,j}^\theta}{\sigma(\gamma)}\right)^2 + (n_{i,j}^\varphi)^2}, \quad \vec{\nu}_{i,j}(\gamma) = \begin{bmatrix} \dfrac{n_{i,j}^\theta}{\sigma(\gamma)\delta_{i,j}(\gamma)} \\ \dfrac{n_{i,j}^\varphi}{\delta_{i,j}(\gamma)} \end{bmatrix}. \quad (32)$$

Finally,

$$D_{i,j}^-(\gamma) = R_{\vec{\nu}_{i,j}(\gamma)}^{-1} \cdot D^-(R_{\vec{\nu}_{i,j}(\gamma)} w_{i,j}^-(\gamma), \eta_{\sigma,i,j}^-; R_{\vec{\nu}_{i,j}(\gamma)} w_{i,j}^+(\gamma), \eta_{\sigma,i,j}^+). \tag{33}$$

Here, $D^-(w_{\sigma,l}, \eta_{\sigma,l}; w_{\sigma,r}, \eta_{\sigma,r})$ is the fluctuation function of a first-order path-conservative method (see [Par06, CMP17]) for the 1d system of balance laws:

$$\partial_t U_\sigma + \partial_x F(U_\sigma) + T^p(h) \partial_x \eta_\sigma = 0, \tag{34}$$

where the rotational invariance-like properties have been used. Notice that (34) is nothing but the 1d shallow water system, with gravity constant $g/\sigma$ and a transport equation for the tangential velocity $\frac{Q_{\vec{p}\perp}}{h_\sigma}$.

**Remark 2.** *The source term $G_\varphi(W)\partial_\varphi \sigma$ does not appear in (34) since $\sigma(\mathbf{x}) = \cos(\varphi)$ is continuous.*

Next, we follow our earlier work[COP17] to define the reconstruction operator $p_{f,i}(\mathbf{x})$ on the cell $\Omega_i$ as

$$p_{f,i}(\mathbf{x}) = p_{f,i}(\mathbf{x}; \{f_j\}_{j \in \mathcal{S}_i}),$$

which is a reconstruction of the data $f_j$, defined by

$$f_j = \eta_{\sigma,j} - \bar{\eta}_i \bar{\sigma}_j, \quad j \in \mathcal{S}_i,$$

where

$$\eta_{\sigma,i} = h_{\sigma,i} - H_{\sigma,i} \quad \text{and} \quad \bar{\eta}_i = \frac{\eta_{\sigma,i}}{\sigma_i}.$$

Note that $f_j$, measure the distance between the cell averages in the stencil of the $i$-th cell and those of the water at rest solution $u_\theta = 0$, $u_\varphi = 0$, $\eta_\sigma = \bar{\eta}_i \cos(\varphi)$. Now, $p_{\eta_\sigma,i}(\mathbf{x})$ is defined in terms of $p_{f,i}(\mathbf{x})$ as follows:

$$p_{\eta_\sigma,i}(\mathbf{x}) = \bar{\eta}_i \cos(\varphi) + p_{f,i}(\mathbf{x}).$$

Despite $\bar{\eta}_i$ being only a second order accurate approximation of the cell average of the free surface in $\Omega_i$, thanks to the continuity of $\sigma = \cos(\varphi)$, it can be easily shown that the accuracy of this modified reconstruction operator for $\eta_{\sigma,i}$ is equal to any standard reconstruction operator and it is trivially well-balanced as $p_{f,i}(\mathbf{x}) = 0$ for any stationary solution (28).

Note that if the shallow-water on Cartesian coordinates is considered, then $p_{\eta_\sigma,i}(\mathbf{x})$ and $p_{f,i}(\mathbf{x})$ differ in an additive constant. Therefore $p_{f,i}(\mathbf{x})$ is replaced by $p_{\eta_\sigma,i}(\mathbf{x})$ in (29).

It is known[COP17] that the resulting method is well-balanced supposing that $D^-(w_{\sigma,l}, \eta_{\sigma,l}; w_{\sigma,r}, \eta_{\sigma,r})$ is a well-balanced first order path-conservative scheme for the standard shallow-water system. In this paper we employ the well-balanced HLLC solver for the standard shallow-water system.

In practice, the integrals appearing in (29) are approximated by quadrature formulas whose order is at least the one of the reconstruction operators:

$$
\begin{aligned}
w'_i(t) &= -\frac{1}{R|\Omega_i|}\Big(\sum_{j\in\mathcal{N}_i}\Big(|E_{i,j}|\sum_{l=0}^{k}\alpha_l F_{\vec{n}_{i,j}}(W^-_{i,j}(\gamma^l_{i,j})) + |E_{i,j}|\sum_{l=0}^{k}\alpha_l \delta_{i,j}(\gamma^l_{i,j})D^-_{i,j}(\gamma^l_{i,j})\Big) \\
&\quad + \sum_{l=0}^{K}\beta_l\left(T^p_\theta(P_i(\mathbf{x}^l_i))\partial_\theta p_{f,i}(\mathbf{x}^l_i) + T^p_\theta(P_i(\mathbf{x}^l_i))\partial_\varphi p_{f,i}(\mathbf{x}^l_i)\right) \\
&\quad + \sum_{l=0}^{K}\beta_l\left(G^1_\varphi(P_i(\mathbf{x}^l_i)) + G^2_\varphi(p_{h_\sigma}(\mathbf{x}^l_i), p_{f,i}(\mathbf{x}^l_i), \sigma(\mathbf{x}^l_i))\right)\partial_\varphi \sigma(\mathbf{x}^l_i)\Big),
\end{aligned}
\qquad(35)
$$

where $\{\gamma^l_{i,j}\}^k_{l=0}$, $\{\alpha_l\}^k_{l=0}$ are the quadrature points and weights of the formula chosen on $E_{i,j}$, and $\{\mathbf{x}^l_i\}^K_{l=0}$, $\{\beta_l\}^K_{l=0}$ those of the formula chosen in $\Omega_i$.

The well-balanced property is preserved: in fact, the integrands vanish at every point on a stationary solution, so that the approximations of the integrals by quadrature formulas also vanish.

The volume quadrature formula may be also used to compute the cell averages of the variable

$$\bar{w}_i = \sum_{l=0}^{K}\beta_l w(\mathbf{x}_l).$$

In that case, the modified reconstruction operator is still well-balanced provided that the averages of $\cos(\varphi)$ are computed by using the same quadrature formula, as it can be easily checked.

The discretization in time is then performed by applying to the ODE system (35) a high order numerical method: TVD RK method will be considered here[GS98]. Finally, the standard CFL condition is considered to obtain a linearly $L^\infty$ stable numerical scheme.

## 4.1 Treatment of wet/dry fronts

Close to a wet-dry transition it is of course important to maintain the positivity of the reconstructed water height. In order to achieve this, the CWENO reconstruction must be slightly modified following [ZS11] and we also follow the ideas described in [GPC07] to modify the 1D HLLC solver. Let us briefly describe the reconstruction procedure in a wet/dry front.

Let us define a cell to be dry if the water height is below some threshold, say $h_\epsilon = 10^{-8}$. First, in a dry cell, the reconstruction is not performed at all and the reconstruction polynomial is defined to be flat.

If a dry cell is present in the stencil $\mathcal{S}_i$ for the reconstruction of a wet cell, the reconstruction order is lowered in the following way. The high order $P_{\mathsf{opt}}$ polynomial is not computed and we set $\hat{d}_0 = 0$ and $P_0(x,y) = 0$. For each low degree candidate polynomial $P_r$ $(r = 1,\ldots,4)$ set $\hat{d}_r = 0$ if the stencil of $P_r$ contains a dry cell and $\hat{d}_r = d_r$ if all cells in the stencil are wet. If at least one of the $\hat{d}_0, \ldots, \hat{d}_4$ coefficient is nonzero, perform the CWENO reconstruction as usual, using the modified linear coefficients, which effectively computes a reconstruction using only polynomials computed from (sub)stencils composed only of wet cells. Otherwise, in case no such all-wet (sub)stencil exists, the reconstruction polynomial is chosen to be flat.

Furthermore, the reconstruction is evaluated at all points needed by the spatial discretization scheme (both inner and boundary quadrature points for the cell). If any of the reconstructed values is negative, the the reconstruction polynomial is replaced by $R(x) \to \overline{u}_j + \theta(R(x) - \overline{u}_j)$ for a suitable choice of $\theta$ that can ensure the positivity of the reconstructed water height. In particular, following [ZS11], letting $h_{\min} = \min_i(h(x_{\xi_i}))$ be the smallest value of the reconstructed water height at all quadrature points $(\xi_i)$, then we set $\theta = \min\left(\left|(h_{\xi_i} - h_\epsilon)/(h_{\xi_i} - h_{\min})\right|, 1\right)$.

## 5 Numerical experiments

In this section we present some numerical tests in order to check the well-balancing and high order properties of the scheme employing the CWENO reconstruction. Here only structured grids on $(x, y)$ or $(\theta, \varphi)$ plane are considered. We use the HLLC scheme written as a PVM method[dlACFN$^+$13]. We use the three step TVD RK method [GS98] that is also third order accurate in time. The constant $d_0$ and $d_1, \ldots, d_4$ appearing in reconstruction are set to 0.75 and 0.0625, respectively. Gauss quadrature formulas are used for the line and volume integrals, using two nodes per direction for the third order scheme and the three nodes per direction for the fourth order schemes. The gravitational constant is set to $g = 9.81 \ m/s^2$. The CFL condition reads as follows in the case of shallow-water on the sphere:

$$\Delta t = \text{CFL} \min_i \left\{ \frac{R \Delta_\theta \Delta_\varphi \cos(\varphi_i)}{(|u_{\theta,i}| + \sqrt{gh_i})\Delta_\varphi + (|u_{\varphi,i}| + \sqrt{gh_i})\Delta_\theta} \right\}, \quad 0 \leq CFL \leq 1, \tag{36}$$

where $\Delta_\theta$ and $\Delta_\varphi$ are the mesh sizes in the $\theta$ and $\varphi$ directions. Here we use CFL = 0.5.

In order to speedup the simulations, a parallel GPU implementation has been performed [GOdlAM11, MdlAC16].

### 5.1 Convergence tests

For this test we consider a vortex perturbing a flat lake surface with water height $h_0 = 2$ over a flat bottom topography where $H = 0$. In particular we consider the initial condition

$$h = h_0 - \frac{\overline{v}^2}{4\alpha g} e^{2.0\alpha(1-r^2)}, \quad u = -\overline{v} e^{\alpha(1.0-r^2)} \frac{ry}{r+\epsilon}, \quad v = \overline{v} e^{\alpha(1.0-r^2)} \frac{rx}{r+\epsilon},$$

|     | water height ||  momentum x  ||  momentum y  ||
| N   | err     | rate | err     | rate | err     | rate |
| --- | ------- | ---- | ------- | ---- | ------- | ---- |
| 25  | 6.8e-1  |      | 1.09e1  |      | 1.09e1  |      |
| 50  | 4.4e-1  | 0.62 | 6.84e0  | 0.68 | 6.84e0  | 0.68 |
| 100 | 1.4e-1  | 1.61 | 1.75e0  | 1.96 | 1.75e0  | 1.96 |
| 200 | 2.55e-2 | 2.50 | 3.2e-1  | 2.45 | 3.2e-1  | 2.45 |
| 400 | 3.41e-3 | 2.91 | 4.16e-2 | 2.94 | 4.16e-2 | 2.94 |

Table 1: Convergence test for the P2/P1 method.

|     | water height ||  momentum x  ||  momentum y  ||
| N   | err     | rate | err     | rate | err     | rate |
| --- | ------- | ---- | ------- | ---- | ------- | ---- |
| 25  | 5.7e-1  |      | 8.26e0  |      | 8.26e0  |      |
| 50  | 2.6e-1  | 1.13 | 3.20e0  | 1.36 | 3.20e0  | 1.36 |
| 100 | 1.81e-2 | 3.83 | 2.6e-1  | 3.60 | 2.6e-1  | 3.60 |
| 200 | 9.56e-4 | 4.24 | 8.06e-3 | 5.03 | 8.05e-3 | 5.03 |
| 400 | 2.42e-5 | 5.30 | 1.97e-4 | 5.35 | 1.97e-5 | 5.35 |

Table 2: Convergence test for the P3/P1 method.

where $r = \sqrt{x^2 + y^2}$ is the distance from the origin and $\alpha = 1$, $\bar{v} = 1$, and $\epsilon = 10^{-16}$. Periodic boundary conditions are set.

The results of the convergence tests are shown in the Tables 1, 2 and 3 for the P2/P1, P3/P1, P3/P2 cases respectively. The P2/P1 scheme achieves the theoretical third order of accuracy at least asymptotically. The two fourth order schemes achieve their theoretical order of convergence already on coarser grids. It is worth noticing that the P3/P2 scheme errors are almost half of the errors of the P3/P1 scheme on all grids with up to 200 × 200 cells. On the other hand, the errors on the grid with 400 × 400 cells are almost identical for the two schemes. This is an indication that on this grid, the nonlinear weights are so close to the optimal linear weights that the reconstruction polynomial is so close to the optimal central third degree polynomial that in the result we cannot distinguish the two schemes. On the other hand, on coarser grids, the reconstruction polynomial is a perturbation of the optimal one by second or third order accurate polynomials and thus the difference of the two schemes show up.

|     | water height ||  momentum x  ||  momentum y  ||
| N   | err     | rate | err     | rate | err     | rate |
| --- | ------- | ---- | ------- | ---- | ------- | ---- |
| 25  | 4.8e-1  |      | 5.98e0  |      | 6.09e0  |      |
| 50  | 1.3e-1  | 1.89 | 1.57e0  | 1.93 | 1.55e0  | 1.98 |
| 100 | 1.00e-2 | 3.73 | 1.2e-1  | 3.77 | 1.2e-1  | 3.74 |
| 200 | 4.88e-4 | 4.32 | 4.67e-3 | 4.62 | 4.88e-4 | 4.63 |
| 400 | 2.42e-5 | 4.33 | 1.97e-4 | 4.57 | 1.97e-5 | 4.57 |

Table 3: Convergence test for the P3/P2 method.

## 5.2 Wet-dry transitions

This numerical test is designed to show the performance of the schemes on solutions where wet/dry fronts appear. Let us consider the paraboloidal topography

$$H(\mathbf{x}) = h_0 \left(1 - \frac{x^2 + y^2}{a^2}\right), \quad \mathbf{x} \in [-2, 2] \times [-2, 2],$$

together with the following periodic analytical solution of the two-dimensional shallow water equations [Tha81]:

$$h(\mathbf{x}, t) = \max\left(0, \frac{\sigma h_0}{a^2}\left(2x \cos(\omega t) + y \sin(\omega t) - \sigma\right) + H(\mathbf{x})\right),$$
$$u_x(\mathbf{x}, t) = -\sigma \omega \sin(\omega t), \quad u_y(\mathbf{x}, t) = \sigma \omega \cos(\omega t),$$

where $u_x$ and $u_y$ are the velocities in the $x$ and $y$ directions, and $\omega = \sqrt{2gh_0}/a$. The values $a = 1$, $\sigma = 0.5$ and $h_0 = 0.1$ have been considered for this test.

The computations have been performed using a quadrilateral mesh with $\Delta x = \Delta y = 0.02$. Comparisons between the numerical and the analytical free surface at points with $y = 0$, $x \in [-2, 2]$ in different times steps are shown in Figure 2. The planar form of the free surface is maintained throughout the computation. Figure 3 shows the comparison between the numerical and the analytical component of the horizontal velocity at the same time steps. Note that only small distortions appear near the shorelines. Moreover, the velocities are maintained after 4 periods. Similar results are obtained for the other component of the velocity.

The results obtained with the P2/P1 and the P3/P1 schemes are very close, both showing only a small oscillation close to the wet dry front, visible mostly in the velocity plot of Figure 3. The results of the P3/P2 scheme have a slightly more pronounced oscillation, that however get smaller with mesh refinement.

## 5.3 Shallow-water on the sphere: water at rest simulation

Here, we consider the shallow-water system on the sphere and the objective of this numerical test is to check the well-balanced property of the numerical scheme. As in [COP17], we consider the rectangular domain $[-180, 180] \times [-89.5, 89.5]$ in the $\tilde{\theta}$-$\tilde{\varphi}$ plane (in degrees) that corresponds to a sphere with two polar caps. A uniform structured grid with $\Delta_{\tilde{\theta}} = \Delta_{\tilde{\varphi}} = 1°$ is considered. The bathymetry is defined as follows: considering the mean bathymetry

$$H_m(\tilde{\theta}, \tilde{\varphi}) = 2.0 - \cos^2\left(\frac{\pi \tilde{\theta}}{60}\right) \cdot \sin^2\left(\frac{\pi \tilde{\varphi}}{60}\right), \tag{37}$$

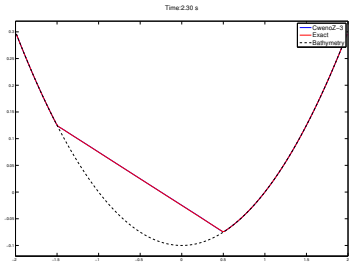
(a) Third order: $t = 2.3$ s

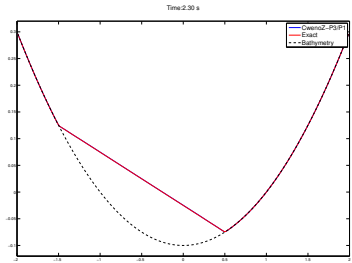
(b) Fourth order: $t = 2.3$ s

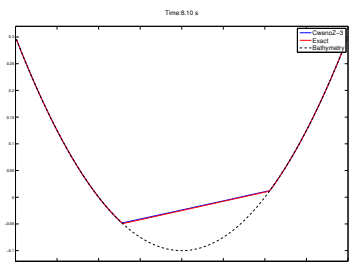
(c) Third order: $t = 8.1$ s

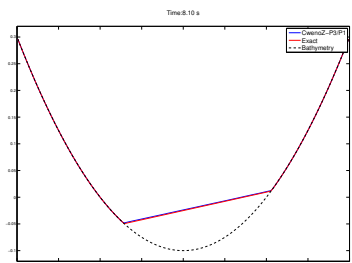
(d) Fourth order: $t = 8.1$ s

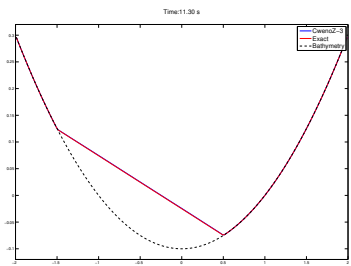
(e) Third order: $t = 11.3$ s

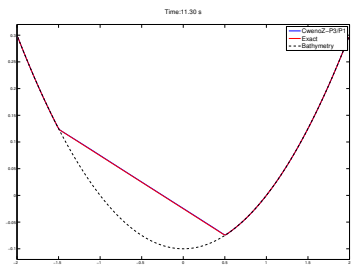
(f) Fourth order: $t = 11.3$ s

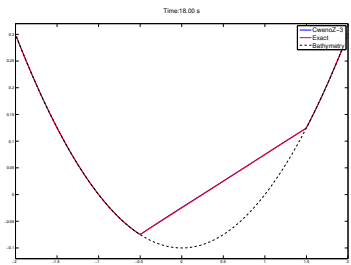
(g) Third order: $t = 18.0$ s

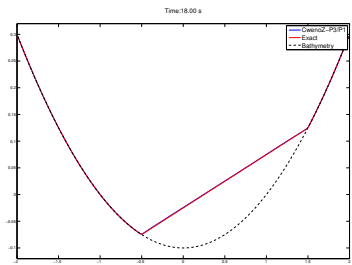
(h) Fourth order: $t = 18.0$ s

Figure 2: 2-d oscillating lake. Surface elevation at $y = 0$, $x \in [-2, 2]$ at different times steps: Exact solution in red, computed solution in blue, bathymetry in doted black line. Third order left column, fourth order (P3/P1) right column.

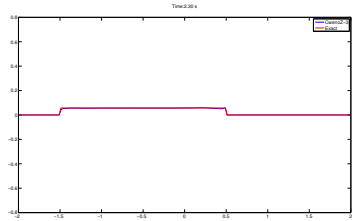
(a) Third order: $t = 2.3$ s

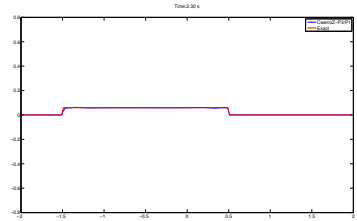
(b) Fourth order: $t = 2.3$ s

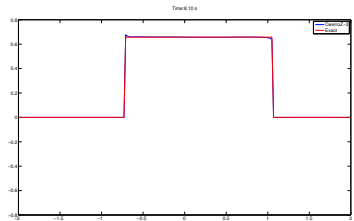
(c) Third order: $t = 8.1$ s

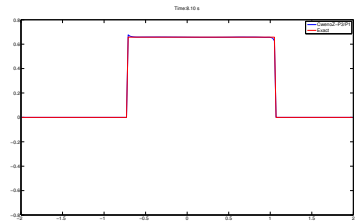
(d) Fourth order: $t = 8.1$ s

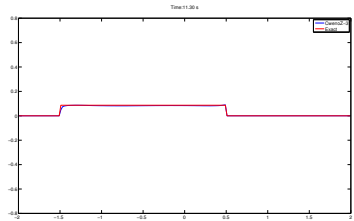
(e) Third order: $t = 11.3$ s

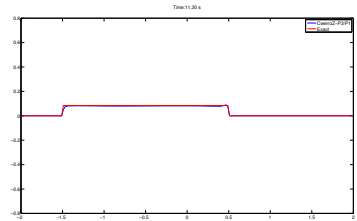
(f) Fourth order: $t = 11.3$ s

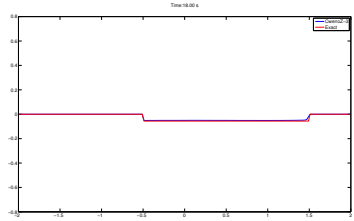
(g) Third order: $t = 18.0$ s

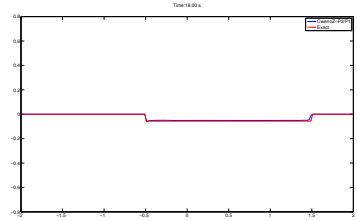
(h) Fourth order: $t = 18.0$ s

Figure 3: 2-d oscillating lake. Horizontal velocity at $y = 0$, $x \in [-2, 2]$ at different times steps: Exact solution in red, computed solution in blue. Third order left column, fourth order (P3/P1) right column.

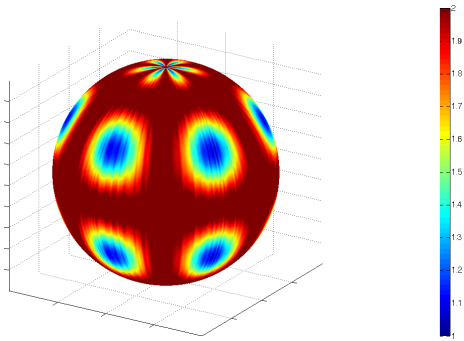

Figure 4: Well-balanced simulation: irregular bathymetry.

| Time | Error (P2/P1) | Error (P3/P1) | Error (P3/P2) |
|------|---------------|---------------|---------------|
| 10m  | 1.47e-15      | 8.63e-16      | 2.57e-15      |
| 60m  | 3.79e-15      | 2.63e-15      | 3.52e-15      |
| 120m | 4.45e-15      | 2.97e-15      | 3.23e-15      |

Table 4: Well-balanced solution: evolution of the error ($L^1$-norm).

we define $H(\tilde{\theta}, \tilde{\varphi})$ by adding to $H_m$ a uniform noise in the interval $[0, 0.2]$. Figure 4 shows the bathymetry for the realization presented here. The initial water depth is set equal to the bathymetry, that is $h(\tilde{\theta}, \tilde{\varphi}, 0) = H(\tilde{\theta}, \tilde{\varphi})$ and $u_{\tilde{\varphi}} = u_{\tilde{\theta}} = 0$. The radius of the sphere is set to $R = 10000$ m. Periodic boundary conditions are prescribed in the eastern and western boundaries, and wall boundary conditions in the northern and southern boundaries corresponding to the polar caps. Table 4 present the evolution of the error in $L^1$-norm for the different high-order numerical schemes. As it can be seen, errors are of the order of the machine accuracy and thus the scheme is exactly well-balanced for the water at rest solution, as expected.

## 5.4 Propagation of a simple wave over an irregular geometry

This test consists on the propagation of a simple wave over an irregular geometry[COP17]. In particular, we consider the same domain but with a finer resolution $\Delta_{\tilde{\theta}} = \Delta_{\tilde{\varphi}} = 0.25°$. The bathymetry $H(\tilde{\theta}, \tilde{\varphi})$ is given by equation (37), the initial water thickness is given by

$$h(\tilde{\theta}, \tilde{\varphi}, 0) = H(\tilde{\theta}, \tilde{\varphi}) + 0.1 e^{-\frac{\tilde{\theta}^2 + \tilde{\varphi}^2}{100}},$$

and initial velocities are set to zero. Again, periodic boundary conditions are prescribed in the eastern and western boundaries, and wall boundary conditions in the northern and southern boundaries. Figure 5 shows the evolution of the free surface computed with the P3/P2 scheme. Note that no spurious oscillations
appear during the wave propagation as the scheme is well-balanced.

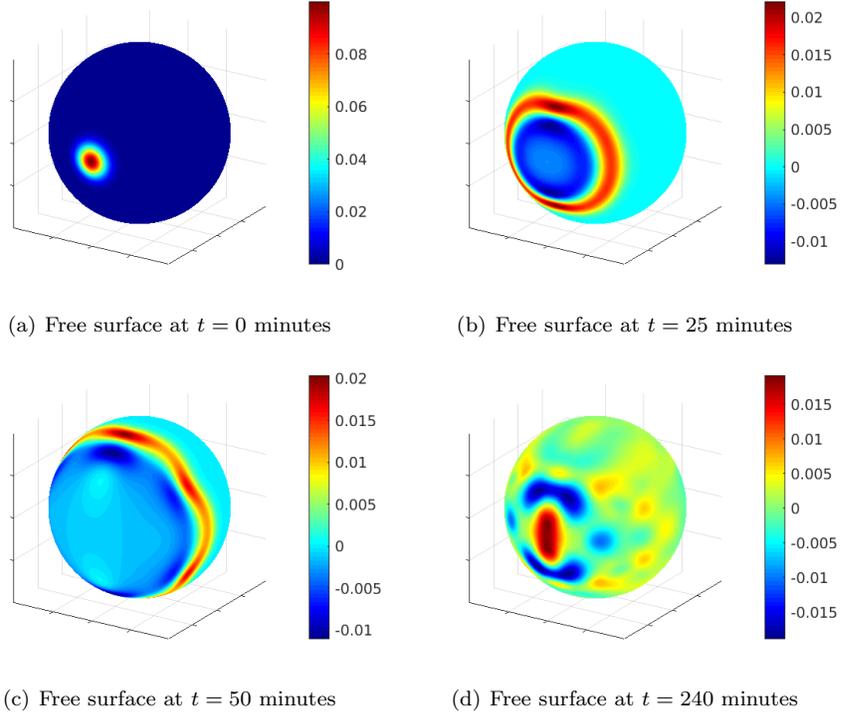

(a) Free surface at $t = 0$ minutes

(b) Free surface at $t = 25$ minutes

(c) Free surface at $t = 50$ minutes

(d) Free surface at $t = 240$ minutes

Figure 5: Evolution of the free surface over an irregular geometry.

In order to compare the third and fourth order schemes, we show in Figure 6 the time series of the water height computed by all three schemes at the point $\theta = 0°, \varphi = 60°$. It can be readily observed that the fourth order schemes compute almost identical solutions, that are indistinguishable in the plot. This is expected, since the solution is smooth and does not have wet-dry fronts. The third order scheme produces a solution that is less accurate: all wave amplitudes are in fact underestimated with respect to the fourth order solutions.

### 5.5 Tohoku event 2011

In this test we simulate the propagation of the 2011 Tohoku tsunami. The main objective of this test is to show that the proposed high order scheme is able to deal effectively with problems set on a real topography. A uniform cartesian grid of the rectangular domain in the $\tilde{\theta}$-$\tilde{\varphi}$ plane given by $[135, 170] \times [25, 49]$ with $\Delta_{\tilde{\theta}} = \Delta_{\tilde{\varphi}} = 4'$ is considered. The mean radius of the Earth is set to $R = 6371009.4$ m and open boundary conditions are prescribed. The topo-bathymetry of the area has been obtained from the ETOPO1 Global Relief

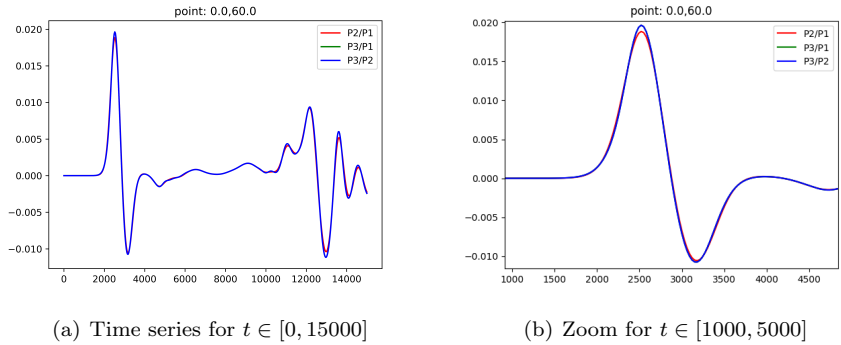

(a) Time series for $t \in [0, 15000]$  (b) Zoom for $t \in [1000, 5000]$

Figure 6: Time series of the free surface at point $(0, 60)$: Comparison with the third and fourth order schemes.

Model. Here we use the initial free surface deformation provided by NCTR-NOAA and the initial velocity field is set to zero. Figures 7 shows the evolution of the tsunami wave propagating near Japan and Figure 8 shows comparison with some DART buoys. Note that the amplitude of the first wave are well captured.

## 6 Conclusion and perspectives

In this paper we presented an high order well-balanced finite volume scheme for the shallow water equations in two space dimensions, both in Cartesian and polar coordinates. The well-balanced property is guaranteed by a path conservative scheme applied to a non-conservative formulation of the equations. High order accuracy in space is achieved by a Central WENO (CWENO) reconstruction procedure from cell averages; time integration was performed with SSP Runge-Kutta methods.

In particular, we presented a third order scheme and two fourth order ones for Cartesian grids. We point out that generalization to even higher accuracy, unstructured and non-conforming meshes is possible, following ideas from [SCR16, DBSR17] to construct appropriate reconstruction procedures.

The scheme were first compared on academic tests to demonstrate the order of accuracy, the non-oscillatory properties in the presence of strong waves and of wet-dry transitions and the well-balanced property with respect to the steady state at rest, both for the Cartesian and for the spherical geometry cases.

The third order reconstruction is based on a central parabola and on four linear polynomials that interpolate data in four directionally-biased substencils. The fourth order generalizations consider a central third order polynomial together with four quadratic (respectively linear) polynomials. It is observed that on regular solutions both fourth order scheme outperform the third order one and that the P3/P2 scheme has lower errors (by a factor of two) than the P3/P1

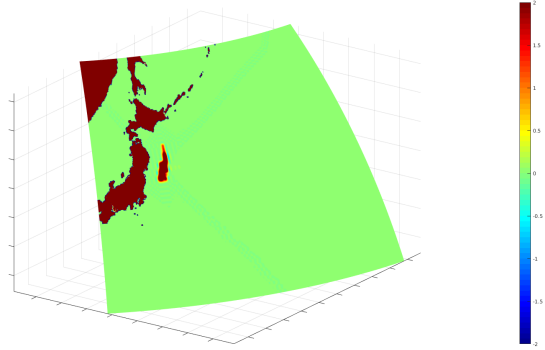

(a) Free surface at $t = 0$ minutes

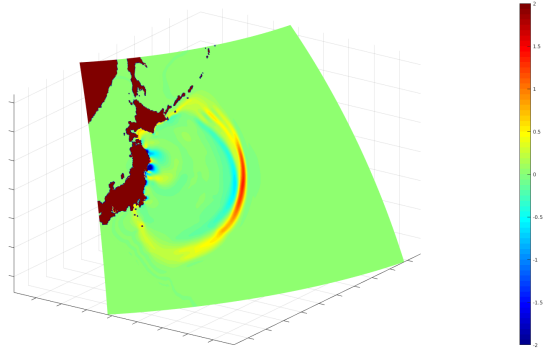

(b) Free surface at $t = 60$ minutes

Figure 7: Free surface evolution for the Tohoku 2011 event.

scheme on coarse meshes, while the two schemes produce almost identical errors on very fine meshes.

Close to discontinuities, all schemes are essentially non-oscillatory, showing small spurious waves and over/undershoots in the velocity field, whose amplitude decrease fast with mesh refinement. More precisely, the two schemes that rely on first degree polynomials (the third order P2/P1 and the fourth order P3/P1) show less pronounced spurious oscillations than the P3/P2 scheme.

Finally, the Tohoku 2011 was simulated, showing the accuracy of the proposed scheme on a problem with a real bathymetry.

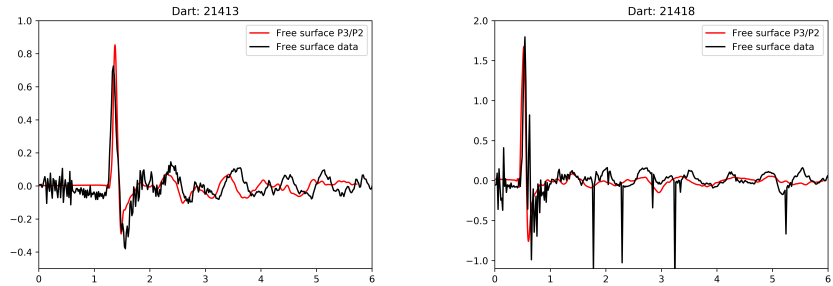

(a) Free surface time evolution at DART buoy # 21413.
(b) Free surface time evolution at DART buoy # 21418.

Figure 8: Comparison with the simulated time series at different DART buoys for the Tohoku 2011 event.

# Acknowledgments

This research has been supported by the Spanish Government and FEDER through the research project MTM2015-70490-C2-1-R. Some of the material used for this work has been discussed during the SHARK-FV 2017 conference.

### Conflict of interest

The authors declare no potential conflict of interests.

[SCR16]  M. Semplice, A. Coco, and G. Russo, *Adaptive mesh refinement for hyperbolic systems based on third-order compact WENO reconstruction*, J. Sci. Comput. **66** (2016), 692–724.

[SHS02]  J. Shi, C. Hu, and C.-W. Shu, *A technique of treating negative weights in WENO schemes*, J. Computat. Phys. **175** (2002), no. 1, 108–127.

[Shu09]  C.-W. Shu, *High order weighted essentially nonoscillatory schemes for convection dominated problems*, SIAM REVIEW **51** (2009), no. 1, 82–126.

[Tha81]  W. C. Thacker, *Some exact solutions to the nonlinear shallow-water wave equations*, J. Fluid Mech. **107** (1981), 499–508.

[TTD11]  P. Tsoutsanis, V. A. Titarev, and D. Drikakis, *WENO schemes on arbitrary mixed-element unstructured meshes in three space dimensions*, J. Computat. Phys. **230** (2011), 1585–1601.

[Zah09]  Y. H. Zahran, *An efficient WENO scheme for solving hyperbolic conservation laws*, Appl. Math. Comput. **212** (2009), 37–50.

[ZQ17]  J. Zhu and J. Qiu, *A new type of finite volume WENO schemes for hyperbolic conservation laws*, J. Sci. Comput. **73** (2017), 1338–1359.

[ZS09]  Y. T. Zhang and C. W. Shu, *Third order WENO scheme on three dimensional tetrahedral meshes*, Comm. Computat. Phys. **5** (2009), 836–848.

[ZS11]  X. Zhang and C.-W. Shu, *Maximum-principle-satisfying and positivity-preserving high-order schemes for conservation laws: survey and new developments*, Proc. R. Soc. Lond. Ser. A Math. Phys. Eng. Sci. **467** (2011), no. 2134, 2752–2776.